\documentclass{amsart}
              [1997/02/02 v1.2e PROC Author Class]

\copyrightinfo{2006}
  {American Mathematical Society}

\newtheorem{theorem}{Theorem}[section]

\newtheorem{proposition}[theorem]{Proposition}

\begin{document}

\title{Compact subspace of products of linearly ordered spaces and co-Namioka spaces}

\author{V.V.Mykhaylyuk}
\address{Department of Mathematics\\
Chernivtsi National University\\ str. Kotsjubyn'skogo 2,
Chernivtsi, 58012 Ukraine}
\email{vmykhaylyuk@ukr.net}



\commby{Ronald A. Fintushel}



\begin{abstract}
It is shown that for any Baire space $X$, linearly ordered compact spaces $Y_1,\dots, Y_n$, compact space $Y\subseteq Y_1\times\cdots \times Y_n$ such that for every parallelepiped $W\subseteq Y_1\times\cdots \times Y_n$ the set $Y\cap W$ is connected, and separately continuous mapping $f:X\times Y\to\mathbb R$ there exists a dense in $X$ $G_\delta$-set $A\subseteq X$ such that $f$ is jointly  continuous at every point of $A\times Y$.
\end{abstract}

\maketitle
\section{Introduction}

Investigation of joint continuity points set of separately continuous function was started by R.~Baire in his classical work \cite{Baire} where he considered functions of two real variables.  The Namioka's result \cite{N} become the impulse to the intensification of these investigations. These result leads to appearance of the following notions which were introduced in \cite{D}.

Let $X$, $Y$ be topological spaces. We say that a separately continuous function $f:X\times Y\to \mathbb R$ {\it has the Namioka property} if there exists a dense in $X$ $G_\delta$-set $A\subseteq X$ such that $f$ is jointly continuous at every point of set $A\times Y$.

A compact space $Y$ is called {\it a co-Namioka space} if for every Baire space $X$ each separately continuous napping $f:X\times Y\to \mathbb R$ has the Namioka property.

Very general results in the direction of study of co-Namioka space properties were obtained in \cite{Bu1, Bu2}. It was obtained in \cite{Bu1, Bu2} that the class of compact co-Namioka spaces is closed over products and contains all Valdivia compacts. Moreover, it was shown in [3] that the linearly ordered compact $[0,1]\times \{0,1\}$ with the lexicographical order is co-Namioka and it was reproved in \cite{Bu1}  (result from \cite{De}) that every completely ordered compact is co-Namioka. These results were generalized in \cite{M}. It was shown in \cite{M} that every linearly ordered compact is co-Namioka.

However, an example given by M..~Talagrand in \cite{T} of a compact space which is not co-Namioka indicates that a closed subspace of co-Namioka compact can be not co-Namioka. Therefore the following question naturally arises: is every compact subspace $Y$ of the product $Y_1\times\cdots \times Y_n$ of finite family of linearly ordered compacts $Y_k$ a co-Namioka?

In this paper we using an approach from \cite{M} we show that under some additional assumptions on $Y$ the formulated above question has the positive answer.

\section{Continuous mappings on compact subspaces of product linearly ordered spaces}

For a mapping $f:X\times Y\to Z$ and a point $(x,y)\in X\times Y$ we put $$f^x(y)=f_y(x)=f(x,y).$$

For a linearly ordered space $X$ and points $a,b\in X$ with $a\leq b$ by $[a,b]$, $[a,b)$, $(a,b]$ and $(a,b)$ we denote the corresponding intervals.

Let $X$ be a topological space. For a mapping $f:X\to\mathbb R$ and a set $A\subseteq X$ by $\omega_f(A)$ we denote the oscillation 
$$\sup\limits_{x,y\in A}|f(x)-f(y)|$$ of the function $f$ on the set $A$. Moreover, for a point $x_0\in X$ by $\omega_f(x_0)$ we denote the oscillation  $$\inf\limits_{U\in{\mathcal U}}\omega_f(U)$$ of the function $f$ at the point $x_0$, where ${\mathcal U}$ is a system of all neighborhoods of $x_0$ in $X$.

\begin{proposition}\label{pr:3.1}
  Let $X_1,\dots, X_n$ be a linearly ordered spaces, $X\subseteq X_1\times\cdots \times X_n$ be a compact space, $\varepsilon >0$ and $f:X\to \mathbb R$ be a continuous mapping. Then there exists an integer $m\in\mathbb N$ such that for every collection $W_1,\dots, W_m$ of parallelepipeds  $$W_k=[a_1^{(k)},b_1^{(k)}]\times\dots\times [a_n^{(k)},b_n^{(k)}]$$ such that $(a_i^{(k)},b_i^{(k)})\cap(a_i^{(j)},b_i^{(j)})=\emptyset$ for every $i\in\{1,\dots n\}$ and distinct $j,k\in \{1,\dots, m\}$, there exists $k_0\in\{1,\dots m\}$ such that $\omega_f(W_{k_0}\cap X)\leq \varepsilon$.
\end{proposition}

\begin{proof} Without loss of the generality we can propose that $X_1,\dots, X_n$ are compacts. Let  $g:X_1\times\cdots \times X_n\to\mathbb R$ be a continuous extension of the mapping $f$. We choose finite covers ${\mathcal U_1},\dots {\mathcal U_n}$ of spaces $X_1,\dots, X_n$ by open intervals such that  $\omega_g(\overline{U}_1\times\dots\times \overline{U}_n)\leq\varepsilon$ for every $U_1\in{\mathcal U_1},\dots U_n\in{\mathcal U_n}$. For every $i\in \{1,\dots, n\}$ we denote by $A_i$ the set of all ending points of intervals $U\in{\mathcal U}_i$ and put $m=|A_1|+|A_2|+\dots +|A_n|+1$. Show that $m$ is the required.

Let $W_1,\dots, W_m$ be a collection of parallelepipeds $$W_k=[a_1^{(k)},b_1^{(k)}]\times\dots\times [a_n^{(k)},b_n^{(k)}]$$ such that $(a_i^{(k)},b_i^{(k)})\cap(a_i^{(j)},b_i^{(j)})=\emptyset$ for every $i\in\{1,\dots n\}$ and distinct $j,k\in \{1,\dots, m\}$. Suppose that $\omega_f(W_{k}\cap X)> \varepsilon$ for every $k\in\{1,\dots m\}$. Then $W_k\not\subseteq \overline{U}_1\times\dots\times \overline{U}_n$ for every $k\in\{1,\dots m\}$ and every $U_1\in{\mathcal U_1},\dots U_n\in{\mathcal U_n}$. Therefore for every $k\in\{1,\dots m\}$ there exist $i\in\{1,\dots n\}$ and $a\in A_i$ such that $a\in (a_i^{(k)},b_i^{(k)})$. Since $m>|A_1|+|A_2|+\dots +|A_n|$, there exist $i\in\{1,\dots n\}$, $a\in A_i$ and distinct $j,k\in \{1,\dots, m\}$ such that $a\in (a_i^{(k)},b_i^{(k)})\cap(a_i^{(j)},b_i^{(j)})$. But it is impossible.
\end{proof}

\section{Main result}

\begin{theorem}\label{th:4.1}
Let $Y_1,\dots, Y_n$ be linearly ordered spaces, $Y\subseteq Y_1\times\cdots \times Y_n$ be a compact subspace such that for every (possibly empty) parallelepiped $$W=[a_1,b_1]\times\dots\times [a_n,b_n]\subseteq Y_1\times\cdots \times Y_n$$ the set $Y\cap W$ is connected. Then $Y$ is co-Namioka.
\end{theorem}

\begin{proof} Let $X$ be a Baire space and $f:X\times Y\to\mathbb R$ be a separately continuous function. We prove that the mapping $f$ has the Namioka property.

We fix $\varepsilon>0$ and show that the open set $$G_{\varepsilon}=\{x\in X: \omega_f(x,y)< \varepsilon
\mbox{\,\,for\,\,every}\,\,y\in Y\}$$ is dense in $X$.

Let $U$ is an open in $X$ nonempty set. Show that there exists an open nonempty set $U_0\subseteq U\cap G_{\varepsilon}$. For every $x\in U$ we denote by $M(x)$ the set of all integers $m\in {\mathbb N}$ for which there exists collection $W_1,\dots, W_m$ of parallelepipeds $$W_k=[a_1^{(k)},b_1^{(k)}]\times\dots\times [a_n^{(k)},b_n^{(k)}]$$ such that $(a_i^{(k)},b_i^{(k)})\cap(a_i^{(j)},b_i^{(j)})=\emptyset$ for every $i\in\{1,\dots n\}$ and distinct $j,k\in \{1,\dots, m\}$ and $\omega_{f^x}(W_{k}\cap Y)> \frac{\varepsilon}{6(n+1)}$ for every $k\in \{1,\dots, m\}$. According to Proposition \ref{pr:3.1} all sets $M(x)$ are upper bounded. For every $x\in U$ we put $\varphi(x)=\max M(x)$, if $M(x)\ne\emptyset$, and $\varphi(x)=0$, if $M(x)=\emptyset$.

It follows from the continuity with respect to the first variable of the function $f$ that for every parallelepiped $W\subseteq Y_1\times\cdots \times Y_n$ the set $\{x\in X:\omega_{f^x}(W\cap Y)>\frac{\varepsilon}{6(n+1)}\}$ is open in $X$. Therefore for every nonnegative $m\in\mathbb Z$ the set $\{x\in U:\varphi(x)>m\}$ is open in $U$, i.e. the function $\varphi:U\to {\mathbb Z}$ is lower semicontinuous on the Baire space $U$. According to \cite{CT} the function $\varphi$ is pointwise discontinuous, i.e. $\varphi$ is continuous at every point from some dense in $U$ sets. There exist an open in $U$ nonempty set $U_1$ and nonegative real $m\in {\mathbb Z}$ such that $\varphi(x)=m$ for every $x\in U_1$.

If $m=0$, then $|f(x,a)-f(x,b)|\leq \frac{\varepsilon}{6(n+1)}<\frac{\varepsilon}{3}$ for every $x\in U_1$ and $a,b\in Y$. Then taking a point $y_0\in Y$ and an open nonempty set  $U_0\subseteq U_1$ such that $\omega_{f_{y_0}}(U_0)<\frac{\varepsilon}{3}$, we obtain that $\omega_f(U_0\times Y)< \varepsilon$. In particular, $U_0\subseteq G_\varepsilon$.

Now we consider the case of $m\in \mathbb N$. Take a point $x_0\in U_1$ and choose a collection $W_1,\dots, W_m$ of parallelepipeds $$W_k=[a_1^{(k)},b_1^{(k)}]\times\dots\times [a_n^{(k)},b_n^{(k)}]$$ such that $(a_i^{(k)},b_i^{(k)})\cap(a_i^{(j)},b_i^{(j)})=\emptyset$ for every  $i\in\{1,\dots n\}$ and distinct $j,k\in \{1,\dots, m\}$ and $\omega_{f^{x_0}}(W_{k}\cap Y)> \frac{\varepsilon}{6(n+1)}$ for every $k\in \{1,\dots, m\}$. Using the continuity of $f$ with respect to the first variable we choose an open neighborhood $U_0\subseteq U_1$ of $x_0$ in $U$ such that $\omega_{f^{x}}(W_{k}\cap Y)> \frac{\varepsilon}{6(n+1)}$ for every $k\in \{1,\dots, m\}$ and $x\in U_0$.

Further without loss of the generality we can propose that $Y_1,\dots, Y_n$ are compacts, i.e. $Y_1=[a_1,b_1],\dots Y_n=[a_n,b_n]$.n For every $i\in\{1,\dots n\}$ we put $$A_i=\{a_i^{(k)}, b_i^{(k)}:1\leq k\leq m\}\cup\{a_i,b_i\}$$ and denote by ${\mathcal V}_i$ the set of all nonempty intervals $[a,b]\subseteq Y_i$ such that $a,b\in A_i$. Moreover, we put $${\mathcal W}={\mathcal V}_1\times\cdots\times {\mathcal V}_n.$$ Show that $\omega_{f^x}(W\cap Y)\leq\frac{\varepsilon}{6}$ for all $x\in U_0$ ³ $W\in {\mathcal W}$.

Suppose that $\omega_{f^x}(W\cap Y)>\frac{\varepsilon}{6}$ for some $x\in U_0$ and $W=[c_1,d_1]\times\cdots\times [c_n,d_n]\in {\mathcal W}$. We choose $z_0=(z_1^{(0)}, \dots,z_n^{(0)}), z_{n+1}=(z_1^{(n+1)}, \dots,z_n^{(n+1)})\in W\cap Y$ such that $|f(x,z_0)-f(x,z_{n+1})|=q>\frac{\varepsilon}{6}$. For certainty we propose that $z_1^{(0)}\leq z_1^{(n+1)},\dots , z_n^{(0)}\leq z_n^{(n+1)}$. Using the continuity of $f^x$ and the fact that for every parallelepiped $V\subseteq Y_1\times\cdots \times Y_n$ the set $V\cap Y$ is connected, we choose points $z_1=(z_1^{(1)}, \dots,z_n^{(1)}),\dots , z_{n}=(z_1^{(n)}, \dots,z_n^{(n)})\in W\cap Y$ such that $z_i^{(0)}\leq z_i^{(1)}\leq z_i^{(2)}\leq\cdots\leq z_i^{(n)}\leq z_i^{(n+1)}$ for every $i\in \{1,\dots, n\}$ and  $|f(x,z_{k-1})-f(x,z_{k})|=\frac{q}{n+1}$ for every $k\in \{1,\dots, n+1\}$. Now for every $i\in \{1,\dots, n\}$ and $k\in \{1,\dots, n+1\}$ we put $c_i^{(k)}=z_i^{(k-1)}$,  $d_i^{(k)}=z_i^{(k)}$ and $V_k=[c_1^{(k)},d_1^{(k)}]\times\dots\times [c_n^{(k)},d_n^{(k)}]$. Note that $(c_i^{(k)},d_i^{(k)})\cap(c_i^{(j)},d_i^{(j)})=\emptyset$ for every $i\in\{1,\dots n\}$ and distinct $j,k\in \{1,\dots, n+1\}$ and $\omega_{f^x}(V_{k}\cap Y)\geq \frac{q}{n+1} > \frac{\varepsilon}{6(n+1)}$ for every $k\in \{1,\dots, n+1\}$.

Since for every $i\in\{1,\dots n\}$ the set $\{k\leq m: (a_i^{(k)},b_i^{(k)})\cap (c_i,d_i)\ne\emptyset\}$ contains at most one element, for the set $$N=\{k\leq m: (a_i^{(k)},b_i^{(k)})\cap (c_i,d_i)=\emptyset\,\,\forall i=1,\dots n\}$$ we have $|N|\geq m-n$. Therefore the system ${\mathcal P}=\{W_k:k\in N\}\cup\{V_j:1\leq j\leq n+1\}$ contains at least $m+1$ parallelepipeds. Besides, $\omega_{f^x}(P\cap Y) > \frac{\varepsilon}{6(n+1)}$ for every $P\in{\mathcal P}$ and $(t_i,u_i)\cap (v_i,w_i)=\emptyset$ for every $i\in \{1,\dots, n\}$, where $$[t_1,u_1]\times\cdots\times [t_n,u_n], [v_1,w_1]\times\cdots\times [v_n,w_n]$$ are distinct parallelepipeds with ${\mathcal P}$. But this contradicts to $\varphi(x)=m$. Thus, $\omega_{f^x}(W\cap Y)\leq\frac{\varepsilon}{6}$ for all $x\in U_0$ and $W\in {\mathcal W}$.

Fix $y\in Y$ and $x\in U_0$. Put ${\mathcal W}_y=\{W\in{\mathcal W}:y\in W\}$ and $V_y=\bigcup\limits_{W\in{\mathcal W}_y}(W\cap Y)$. Clearly that $V_y$ is a neighborhood of $y$ in $Y$. Since $\omega_{f^x}(W\cap Y)\leq\frac{\varepsilon}{6}$ for every $W\in {\mathcal W}_y$ and $y\in \bigcap\limits_{W\in{\mathcal W}_y}(W\cap Y)$, $\omega_{f^x}(V_y)\leq \frac{\varepsilon}{3}$. Now using the continuity of $f$ with respect to the first variable at $(x,y)$ we choose a neighborhood $\tilde{U}$ of $x$ in $X$ such that $\omega_{f_{y}}(\tilde{U})<\frac{\varepsilon}{3}$. Then $\omega_f(\tilde{U}\times V_y)< \varepsilon$, in particular, $\omega_f(x,y)< \varepsilon$ for every $y\in Y$ and $x\in U_0$. Thus, $U_0\subseteq G_\varepsilon$.

\end{proof}

\bibliographystyle{amsplain}

\end{document}